\documentclass[12pt]{amsart}

\usepackage{amssymb,amsthm,amsmath,scalefnt}
\usepackage[numbers,sort&compress]{natbib}
\usepackage{color}
\usepackage{graphicx}
\usepackage{tikz}
\usepackage{amssymb,amsthm,amsmath}
\usepackage[numbers,sort&compress]{natbib}
\usepackage{color}
\usepackage{graphicx}
\usepackage{tikz}
\usepackage{xparse}
\usepackage{url}
\usepackage[colorlinks]{hyperref}
\usepackage{appendix}
\usepackage{listings}
\definecolor{maroon}{RGB}{133, 5, 63}
\definecolor{teal}{RGB}{0, 128, 96}
\definecolor{forestgreen}{RGB}{34, 139, 34}
\lstdefinelanguage{julia}{
basicstyle=\small\ttfamily,
alsoletter=",
classoffset=1,
comment=[l]{\#},
commentstyle=\color{gray},
keywords={solve, differentiate, subs, sub, real_solutions, det, sum, max, conditional_count, map, filter, map_filter, union, zip, product, is_finite, is_successful, is_nonsingular, track, is_real, is_success, first, evaluate, flatten, bitmask_filter, accumulate, polynomial_interpolants, reverse, coefficients, vcat, solution_from_necklace, length, prod, compress, System, certify total_degree_start_solutions, degrees, variables, iterate, struct, collect, convert, stretched_cube, stretched_cubes, weight_vector, weight_vectors, perm_to_segments, perm_to_mixedcell, mixed_cell_iterator, rand_approx_unit, norm, fixed, polyhedral_system, reduce, findall, println, eachrow, eval, parse, readlines, certify, certificates},
keywordstyle={\color{teal}},
classoffset=2,
morekeywords={@var, @time, for, end, if, while, else, begin, global, in },
keywordstyle={\color{maroon}},
classoffset=3,
morekeywords={using, function, return, const},
keywordstyle={\color{blue}},
classoffset=4,
morekeywords={julia, >},
keywordstyle={\color{forestgreen}},
xleftmargin=0.2cm,
xrightmargin=1em,
columns=fullflexible,
keepspaces=true,
}
\lstdefinelanguage{M2}{
basicstyle=\small\ttfamily,
alsoletter=",
classoffset=1,
comment=[l]{\--},
commentstyle=\color{gray},
keywords={solve, differentiate, subs, sub, real_solutions, det, sum, max, count, conditional_count, map, filter, map_filter, union, zip, product, is_finite, is_successful, is_nonsingular, track, is_real, is_success, first, evaluate, flatten, bitmask_filter, entries, unique, join, accumulate, polynomial_interpolants, reverse, coefficients, vcat, solution_from_necklace, length, prod, compress, total_degree_start_solutions, degrees, variables, iterate, struct, collect, convert, stretched_cube, stretched_cubes, weight_vector, weight_vectors, perm_to_segments, perm_to_mixedcell, mixed_cell_iterator, rand_approx_unit, norm, fixed, polyhedral_system, reduce, findall, eachrow, matrix,  peek, coefficientRing, numgens, ring, print, ideal, toString  },
keywordstyle={\color{cyan}},
classoffset=2,
morekeywords={@var, time, for, end, if, in, from, to,  while, else, begin,list ,new, do, scan, apply, close, openOut},
keywordstyle={\color{maroon}},
classoffset=3,
morekeywords={using, function, return, const, needsPackage},
keywordstyle={\color{blue}},
classoffset=4,
morekeywords={julia, >},
keywordstyle={\color{forestgreen}},
classoffset=5,
morekeywords={ polySystem, point, newton, certifyRealSolution},
keywordstyle={\color{purple}},
xleftmargin=0.2cm,
xrightmargin=1em,
columns=fullflexible,
keepspaces=true,
}

\title[ ]{A Counterexample to Fermi Isospectral Rigidity for Two Dimensional Discrete Periodic Schr\"odinger Operators}
\author{Taylor Brysiewicz}
\address[T. Brysiewicz]{ Department of Mathematics, Western University, London, ON, CA} \email{tbrysiew@uwo.ca}

\author{Matthew Faust}
\address[M. Faust]{ Department of Mathematics, Michigan State University, East Lansing, MI, 48840, USA} \email{mfaust@msu.edu}

\author{Wencai Liu}
\address[W. Liu]{ Department of Mathematics, Texas A\&M University, College Station, TX 77843-3368, USA} \email{ wencail@tamu.edu}



\theoremstyle{plain}
\newtheorem{theorem}{Theorem}[section]

\newtheorem{lemma}[theorem]{Lemma}

\newtheorem{remark}{Remark}
\newcommand{\C}{\mathbb{C}}

\newcommand{\Z}{\mathbb{Z}}
\newcommand{\ZZ}{\mathbb{Z}}

\theoremstyle{definition}
\newtheorem{definition}{Definition}
\newtheorem{question}{Question}
\newtheorem{conjecture}{Conjecture}

\begin{document}
	\begin{abstract}
    Using numerical certification, we prove the existence of  a nontrivial real-valued two dimensional periodic potential  whose associated discrete Schr\"odinger operator is Fermi isospectral to the zero potential.
This provides a negative answer to a question posed by the third author concerning the rigidity of Fermi isospectrality in dimension two. This example also disproves a conjecture of Gieseker, Kn\"orrer, and Trubowitz in the 1990s stating that for any nontrivial real-valued periodic potential in dimension two, the Fermi variety is irreducible at all energy levels.
\end{abstract}
	
	\maketitle

 \section{Introduction and main results}

Let $q=(q_1,\ldots,q_d)\in \Z_+^d$ be pairwise coprime and define the lattice
\[
\Gamma=q_1\Z\oplus q_2 \Z\oplus\cdots\oplus q_d\Z.
\]
A function $V\colon\Z^d\to \C$ is called $\Gamma$-periodic (or simply, periodic) if
\[
V(n+\gamma)=V(n), \qquad \text{ for all }\gamma\in \Gamma,\; n\in\Z^d.
\]
Throughout, we assume the potential $V$ is $\Gamma$-periodic.
For $n\in\Z^d$, denote
\[
\|n\|_1=\sum_{j=1}^d |n_j| \quad \text{ where } \quad  n=(n_1,n_2,\cdots,n_d)
\]
and let $\Delta$ be the discrete Laplacian on the lattice $\Z^d$, defined by
\[
(\Delta u)(n)=\sum_{n'\in\Z^d,\; \|n'-n\|_1=1}u(n').
\]

In this article, we study the isospectrality problem and the irreducibility of Fermi varieties for discrete periodic Schr\"odinger operators of the form $\Delta+V$. 
We refer the reader to the survey articles~\cite{ksurvey,lsur} for background and recent developments on these topics.

\begin{definition}
The {\it Bloch variety} $B(V)$ of $\Delta+V$ consists of all pairs $(k,\lambda)\hspace{-1pt}\in\hspace{-1pt}\C^{d+1}$ for which there exists a nonzero solution of
\begin{equation}
(\Delta u)(n)+V(n)u(n)=\lambda u(n), \qquad n\in\Z^d,
\label{spect_0}
\end{equation}
satisfying the Floquet--Bloch boundary conditions
\begin{equation}
u(n+q_j\mathbf e_j)=e^{2\pi i k_j}u(n), 
\qquad j=1,2,\ldots,d,\; n\in\Z^d,
\label{Fl}
\end{equation}
where $k=(k_1,k_2,\ldots,k_d)\in\C^d$ and $\{\mathbf e_j\}_{j=1}^d$ is  the standard basis of $\Z^d$.
For  fixed $\lambda\in\C$, the {\it Fermi surface} (or {\it Fermi variety}) is defined by
\[
F_{\lambda}(V)=\{k:(k,\lambda)\in B(V)\}.
\]
\end{definition}
Our primary interest, in this paper, concerns the isospectrality of periodic Schr\"odinger operators. 
Various notions of isospectrality have been studied in the literature; see for instance \cite{liu2021fermi,kang} for a discussion of   such problems. 
Among them, Fermi isospectrality appears to be one of the most difficult and least understood. 
This notion, introduced in \cite{liu2021fermi}, is the focus of the present work.

\begin{definition}\footnote{After introducing the Floquet matrices in the next section, one can equivalently define Fermi isospectrality by the identity 
$\det (\mathcal D_X(z)-\lambda_0 I)= \det (\mathcal D_Y(z)-\lambda_0 I)$.}\label{fermiiso} 
\cite{liu2021fermi}
Let $X$ and $Y$ be two $\Gamma$-periodic functions. 
We say that $X$ and $Y$ are {\it Fermi isospectral at the energy level $\lambda_0\in\C$} if
\[
F_{\lambda_0}(X)=F_{\lambda_0}(Y).
\]
They are {\it Fermi isospectral} if they are Fermi isospectral at some~$\lambda_0~\in~\C$.
\end{definition}

In \cite{liu2021fermi}, several rigidity results for discrete periodic Schr\"odinger operators were obtained. 
One of them is the following. 

\begin{theorem}\label{thm}
\cite{liu2021fermi}
Let $d\ge 3$. The only real potential $V$ which is Fermi isospectral to the zero potential $\textbf{0}$ is $\textbf{0}$ itself. 
\end{theorem}

\noindent It is natural to ask whether the same rigidity result holds in dimension $d=2$.

\medskip

\noindent
\begin{question}{~\cite{lsur}.}~\label{question1}
Let $d=2$ and suppose that $V$ is a real periodic potential Fermi isospectral to $\textbf{0}$.  Must it follow that $V = \bf 0$?
\end{question}
\medskip

\noindent In the present paper we give a negative answer to this question.

\begin{theorem}\label{thm1}
Let $d=2$. 
There exists a nontrivial  $3\Z\times5\Z$-periodic  real potential $V$ which is Fermi isospectral to the zero potential at $\lambda=0$, that is,
\[
F_0(V)=F_0(\mathbf 0).
\]
\end{theorem}

Our second goal in this paper concerns the irreducibility of Fermi varieties. 
The irreducibility of Fermi varieties (as well as Bloch varieties) and related applications, such as embedded eigenvalues and spectral band edges, have been extensively studied over the past three decades; see for example
\cite{im,aim,lmt,ktcmh90,bktcm91,bat1,ls,shi2,flscmp,kv06cmp,kvcpde20,shi1,faust,fg25,dksjmp20,LiuQE,flm,mastafa-theo,batcmh92,flm1}.

The following irreducibility result reveals a connection between the irreducibility of Fermi varieties and Fermi isospectrality.
Denote by $[V]$ the average of $V$ over one periodicity cell.

\begin{theorem}\label{gcf1}
\cite{liu1}
For any $d\ge3$, the Fermi variety $F_{\lambda}(V)/\Z^d$ is irreducible for all $\lambda\in\C$.
For $d=2$, the Fermi variety $F_{\lambda}(V)/\Z^2$ is irreducible for all $\lambda\in\C$ except possibly for $\lambda=[V]$, and $F_{[V]}(V)/\Z^2$ has at most two irreducible components. 
Moreover, if the Fermi variety $F_{[V]}(V)/\Z^2$ has two irreducible components, then
\begin{equation}\label{gd2}
  F_{[V]}(V)=F_0(\mathbf 0).  
\end{equation}
In particular, $V$ is Fermi isospectral to the constant potential (equivalently, to the zero potential after shifting).
\end{theorem}

\begin{remark}\label{rem}
    When $d=2$ and $V= \bf 0$, direct computation (see for example \cite{lsur}) shows that $F_0(\bf 0)/\Z^2$ has exactly two irreducible components.
\end{remark}
We note that Theorem \ref{gcf1}  applies to complex-valued potentials $V$. 
It is not difficult to see that there exist nonconstant complex-valued potentials on $\mathbb{Z}^2$ for which \eqref{gd2} holds. 
By Remark \ref{rem}, such potentials imply that $F_{[V]}(V)$ has two irreducible components.
Historically, for real-valued potentials it has been widely believed that constant potentials are the only case for which the Fermi variety $F_{\lambda}(V)/\Z^2$ can be reducible at some energy level. 
This belief was formulated as a conjecture by Gieseker, Kn\"orrer, and Trubowitz in the early 1990s~\cite{GKTBook}.

\medskip

\noindent
\begin{conjecture}{~\cite[p.~43]{GKTBook}.}~\label{conj1}
Let $q_1$ and $q_2$ be two distinct odd primes.
Assume that $V$ is a nonconstant real-valued $q_1\Z\oplus q_2\Z$-periodic potential. 
Then the Fermi variety $F_{\lambda}(V)/\Z^2$ is irreducible for every $\lambda\in\C$.
\end{conjecture}
\medskip
	
Indeed, Conjecture~\ref{conj1} is true if one restricts to ``separable'' potentials; see \cite{Liu_2022}. However, a corollary of Theorem~\ref{thm1} is that Conjecture~\ref{conj1} is, in general, false.

\begin{theorem}\label{thm3}
There exists a nonconstant real-valued periodic potential $V$ on $\Z^2$ such that the Fermi variety $F_{\lambda}(V)/\Z^2$ is reducible for some $\lambda\in\C$.
\end{theorem}

We prove Theorem~\ref{thm1} using standard certification techniques, implemented in  \texttt{Macaulay2}~\cite{M2} and \texttt{julia} \cite{julia}. Our code is given in the Appendix. First we formulate the space of $3\ZZ \times 5\ZZ$-periodic potentials Fermi isospectral to the potential $\textbf{0}$ at energy level $0$ as a solution set to a collection of polynomial equations.  We then  \textit{certify} that a smooth isolated real solution to this polynomial system exists  using the implementation of Krawczyk's method~\cite{Krawczyk} in the \texttt{Macaulay2}~\cite{M2} package  \texttt{NumericalCertification}~\cite{M2Certification}.
By the implicit function theorem, this guarantees the existence of a real curve of potentials Fermi isospectral to $\textbf{0}$ at energy level~$0$. We independently verify this certification with a \texttt{julia}~\cite{julia} implementation of Krawczyk's method given in \texttt{HomotopyContinuation.jl} \cite{HC,JuliaCertification}. 

Our \texttt{Macaulay2} code \emph{proves} that such a real curve exists, however,  it does not indicate how the candidate solution used in the certification process was obtained. We briefly describe how such an approximate solution was found using \texttt{homotopy iterators} \cite{HomotopyIterators}, implemented in \texttt{HomotopyContinuation.jl},  which allow one to bypass the cost of fully numerically solving the system.

\section{Floquet matrices and polynomial formulation}

\subsection{Floquet matrices}

In this section we recall the polynomial formulation of the Floquet problem. 
Our goal is to write, explicitly, the Floquet matrix corresponding to a $3\Z\times5\Z$-periodic potential.

Let $\C^{\ast}=\C\setminus\{0\}$ and write $z=(z_1,z_2,\ldots,z_d)$. 
For any $z\in(\C^{\ast})^d$, consider the equation
\begin{equation}\label{gei1}
(\Delta+V)u=\lambda u
\end{equation}
with boundary conditions
\begin{equation}\label{gei2}
u(n+q_j\mathbf e_j)=z_j\,u(n), 
\qquad j=1,2,\ldots,d,\quad n\in\Z^d .
\end{equation}

We fix a fundamental domain $W$ of $\Gamma$:
\[
W=\{n=(n_1,n_2,\ldots,n_d)\in\Z^d:0\le n_j\le q_j-1,\; j=1,2,\ldots,d\},
\]
and restrict the operator $\Delta+V$ to the $Q=q_1q_2\cdots q_d$-dimensional space
\[
\{u(n):n\in W\}.
\]
Equation \eqref{gei1} with boundary conditions \eqref{gei2} (equivalently \eqref{spect_0} and \eqref{Fl}),
then becomes an eigenvalue problem for a $Q\times Q$ some matrix $\mathcal D_V(z)$ 
(or $D_V(k)$ in the $k$-variable). 
We call $\mathcal D_V(z)$ (respectively $D_V(k)$) the \emph{Floquet matrix}.

We define 
\begin{equation}\label{g16}
\mathcal P_V(z,\lambda)=\det(\mathcal D_V(z)-\lambda I), 
\qquad 
P_V(k,\lambda)=\det(D_V(k)-\lambda I).
\end{equation}
so that the Fermi variety can be written as
\begin{equation}\label{g11}
F_{\lambda}(V)=\{k\in\C^d:P_V(k,\lambda)=0\}.
\end{equation}

The following lemma characterizes Fermi isospectrality.

\begin{lemma}\cite[Lemma 2.3 and Theorem 4.1]{liu2021fermi}\label{keylem}
Two $\Gamma$-periodic functions $X$ and $Y$ are Fermi isospectral at the energy level $\lambda_0\in\C$ if and only if
\[
\det (\mathcal D_X(z)-\lambda_0 I) = 
\det (\mathcal D_Y(z)-\lambda_0 I).
\]
In particular, a periodic potential $V$ is Fermi isospectral to the zero potential at $\lambda_0=0$ if and only if
\begin{equation}\label{gmain}
\det \mathcal D_V(z) =  \det \mathcal D_{\mathbf 0}(z).
\end{equation}
If two $\Gamma$-periodic functions $X$ and $Y$ are Fermi isospectral, then
\begin{equation}\label{gv}
   [X]=[Y] .
\end{equation}
\end{lemma}

Specializing to dimension $d=2$, the Floquet matrix $\mathcal D_V(z)$ has the following expression when $q_1\ge3$ and $q_2\ge3$ (see \cite{GKTBook,liu1}):

\[
(\mathcal{D}_V(z_1,z_2))_{(m,n),(m',n')}=
\begin{cases}
V(m,n) & m=m',\ n=n',\\
1 & (m-m')^2+(n-n')^2=1,\\
z_1 & m'=1,\ m=q_1,\ n=n',\\
z_1^{-1} & m'=q_1,\ m=1,\ n=n',\\
z_2 & m=m',\ n'=1,\ n=q_2,\\
z_2^{-1} & m=m',\ n'=q_2,\ n=1,\\
0 & \text{otherwise}.
\end{cases}
\]

For $(q_1,q_2) = (3,5)$, and the ordering
\[
(1,1),(2,1),(3,1),(1,2),(2,2),(3,2),\dots,(1,5),(2,5),(3,5).
\]
 the  Floquet matrix $\mathcal D_V(z)$ takes an explicit form: writing $v_{i,j}$ for the value $V(i,j)$, we have   that 
\[
\scalefont{1.3}{
\mathcal D_{V}(z) = \left(
\begin{smallmatrix}
v_{1,1} & 1 & z_1^{-1} & 1 & 0 & 0 & 0 & 0 & 0 & 0 & 0 & 0 & z_2^{-1} & 0 & 0 \\
1 & v_{2,1} & 1 & 0 & 1 & 0 & 0 & 0 & 0 & 0 & 0 & 0 & 0 & z_2^{-1} & 0 \\
z_1 & 1 & v_{3,1} & 0 & 0 & 1 & 0 & 0 & 0 & 0 & 0 & 0 & 0 & 0 & z_2^{-1} \\
1 & 0 & 0 & v_{1,2} & 1 & z_1^{-1} & 1 & 0 & 0 & 0 & 0 & 0 & 0 & 0 & 0 \\
0 & 1 & 0 & 1 & v_{2,2} & 1 & 0 & 1 & 0 & 0 & 0 & 0 & 0 & 0 & 0 \\
0 & 0 & 1 & z_1 & 1 & v_{3,2} & 0 & 0 & 1 & 0 & 0 & 0 & 0 & 0 & 0 \\
0 & 0 & 0 & 1 & 0 & 0 & v_{1,3} & 1 & z_1^{-1} & 1 & 0 & 0 & 0 & 0 & 0 \\
0 & 0 & 0 & 0 & 1 & 0 & 1 & v_{2,3} & 1 & 0 & 1 & 0 & 0 & 0 & 0 \\
0 & 0 & 0 & 0 & 0 & 1 & z_1 & 1 & v_{3,3} & 0 & 0 & 1 & 0 & 0 & 0 \\
0 & 0 & 0 & 0 & 0 & 0 & 1 & 0 & 0 & v_{1,4} & 1 & z_1^{-1} & 1 & 0 & 0 \\
0 & 0 & 0 & 0 & 0 & 0 & 0 & 1 & 0 & 1 & v_{2,4} & 1 & 0 & 1 & 0 \\
0 & 0 & 0 & 0 & 0 & 0 & 0 & 0 & 1 & z_1 & 1 & v_{3,4} & 0 & 0 & 1 \\
z_2 & 0 & 0 & 0 & 0 & 0 & 0 & 0 & 0 & 1 & 0 & 0 & v_{1,5} & 1 & z_1^{-1} \\
0 & z_2 & 0 & 0 & 0 & 0 & 0 & 0 & 0 & 0 & 1 & 0 & 1 & v_{2,5} & 1 \\
0 & 0 & z_2 & 0 & 0 & 0 & 0 & 0 & 0 & 0 & 0 & 1 & z_1 & 1 & v_{3,5}
\end{smallmatrix}
\right)}.
\]
\subsection{Reduction to an algebraic system}
Treating the values $V(m,n)$ of the potential as variables $\textbf{v}=(v_{m,n})_{(m,n)=(1,1)}^{(3,5)}$, we construct polynomial equations \[F:f_1(\textbf{v})=\ldots=f_{14}(\textbf{v})=0 \quad \quad \text{ where } \quad \textbf{v}=(v_{1,1},\ldots,v_{3,5})\] obtained by expanding the identity \eqref{gmain} and comparing coefficients. In the variables $z_1$ and $z_2$ there are a total of $29$ monomials, however, several coefficients are identical leaving  only $14$ unique coefficients. These $14$ polynomials $f_1,\ldots,f_{14} \in \mathbb{Q}[v_{1,1},\ldots,v_{3,5}]$ comprise the polynomial system $F$. Code for constructing this system is given in the Appendix.
By Lemma~\ref{keylem}, the problem of Fermi isospectrality reduces to finding a nonzero real solution to $F=\textbf{0}$.

\begin{theorem}\label{thm2}
There exists a nontrivial $3\Z\times5\Z$-periodic real potential $V$ that is Fermi isospectral to the zero potential at energy $\lambda=0$ if and only if there exists a real solution to the system $F=\textbf{0}$ other than $\textbf{0}$.
\end{theorem}
	 
\section{Computational proof of Theorem \ref{thm1}}
The $3\mathbb{Z} \times 5\mathbb{Z}$-periodic potential identified by  matrix
\[
V^*=\begin{pmatrix}
3.587996 & -0.130207 & -1.455985 & -4.537496 & -0.224238 \\
0.688112 & -0.678205 & -0.377625 & 4.597172 & 1.810359 \\
1.870801 & 0.006187 & -4.264991 & -0.362522 & -0.529361
\end{pmatrix}
\] is provably near a smooth real solution $\widetilde{V}$ to the system $F=(f_1,\ldots,f_{14})=\textbf{0}$ of Theorem \ref{thm2}. As a first check, one may compute 
\begin{align*}
\textrm{det}(\mathcal D_{V^*}(z))  &= z_1^5 \hspace{-1pt} + \hspace{-1pt}  z_1^{-5} \hspace{-1pt}  + \hspace{-1pt}  z_2^3  \hspace{-1pt} +  \hspace{-1pt} z_2^{-3} \hspace{-1pt} 
+  \hspace{-1pt} 2.3\times 10^{-8}\Big(
 \hspace{-3pt} - z_1^2
+\cdots 
+5.8\times 10^{-7}z_1^{-3}z_2^{-1}
\Big) \\
&\approx z_1^5 \hspace{-1pt} + \hspace{-1pt}  z_1^{-5} \hspace{-1pt}  + \hspace{-1pt}  z_2^3  \hspace{-1pt} +  \hspace{-1pt} z_2^{-3}  = \det(\mathcal D_{\textbf{0}}(z)).
\end{align*}
However, the above calculation is clearly \textit{not} a proof, due to its numerical nature. Rather, to prove our claims, we rely on a numerical certification procedure known as Krawczyk's method \cite{Krawczyk}. 

We begin by showing how to compute the polynomials $F=(f_1,\ldots,f_{14})$ and find a candidate approximate real solution to $F=\textbf{0}$, like $V^*$ above. Next, we summarize Krawczyk's method and explain how its success proves the existence of a true smooth real solution $\widetilde{V}$ to $F=\textbf{0}$ near $V^*$. Finally, we apply Krawczyk's method to $V^*$ and $F$ using two independent implementations, one in  \texttt{Macaulay2} \cite{M2,M2Certification} and another in \texttt{julia} \cite{julia,HC,JuliaCertification}. 

\subsection{Generating the system} Although fairly direct, one must apply some care in constructing the system $F=\textbf{0}$ of Theorem \ref{thm2}. Na\"ively taking the determinant of $\mathcal D_V(z)$ in \texttt{Macaulay2} is not sufficient and does not terminate in a reasonable amount of time. To make the computation feasible, first, we change the default determinant computation strategy of \texttt{Bareiss} to \texttt{Cofactor} to take advantage of the sparsity of $\mathcal D_{V}(z)$. Secondly, we treat $z_1^{-1}$ and $z_2^{-1}$ as their own formal variables until we consolidate the coefficients of monomials in $\mathcal D_{V}(z) - \mathcal D_{\textbf{0}}(z)$ by working   modulo the ideal $\langle z_1z_1^{-1}-1,z_2z_2^{-1}-1\rangle$. With these adjustments, \texttt{Macaulay2} computes the polynomials $f_1,\ldots,f_{14}$ in less than five seconds. Our \texttt{Macaulay2} code to compute these polynomials and print them to files (in both \texttt{Macaulay2} and \texttt{julia} syntax) can be found in the Appendix.

\subsection*{Finding a candidate solution}
In this subsection, we discuss one way to produce a candidate 
solution to $F=\textbf{0}$. The polynomials $f_1,\ldots,f_{14}$ have degrees 
\[
3, 1, 6, 4, 2, 9, 7, 5, 4, 12, 10, 2, 7, 15
\]
respectively.  B\'ezout's theorem \cite[Section 5]{Fulton1969} then states that $X_{\mathbb{C}}$ has a degree no larger than the product $\prod_{i=1}^{14} \textrm{deg}(f_i) = 4,572,288,000$. Indeed, this is the number of paths one needs to track to solve for a hyperplane slice of $X_{\mathbb{C}}$ via the standard numerical algebraic geometry technique of a \textit{total degree homotopy} \cite{NAG}.  At the rate of approximately $\approx 350$ paths per minute, solving for all solutions would require around $\approx 25$ years. Rather, we use the lazy-evaluation approach of \cite{HomotopyIterators} and construct a \textit{homotopy iterator} for the solution set to $F=\textbf{0}$ (intersected with a random hyperplane). This data structure allows us to track one path at a time until we find a real solution. This process terminated after approximately $400$ minutes, and found a real solution on the $142295^{th}$ path.

\subsection{Krawczyk's method} The key algorithm we use to prove that $F$ has a real nontrivial root is \textit{Krawczyk's method} \cite{Krawczyk}. Fix a square polynomial system 
\[
F\colon f_1,\ldots,f_n \in \mathbb{Q}[x_1,\ldots,x_n]
\]
and a candidate floating point approximate solution $\textbf{x}^* \in \mathbb{C}^n$, Krawczyk's method constructs and tries to prove that a bounding box \[B = [a_1,a_1'] \times [b_1,b_1'] \times \cdots \times [a_n,a_n'] \times [b_n,b_n'] \subseteq \mathbb{R}^{2n} \cong \mathbb{C}^n
\]
containing $\textbf{x}^*$ contracts under the Newton operator $N_F:\mathbb{C}^n \to \mathbb{C}^n$ associated to $F=(f_1,\ldots,f_n)$, i.e. $N_F(B) \subseteq \textrm{interior}(B)$. By \textit{Banach's fixed point theorem} this proves there must be a \textit{unique} fixed point $\widetilde{\textbf{x}}$ of $N_F$ in $B$. Such fixed points are precisely the roots of $F$, and so finding a box which contracts proves the existence of a unique isolated root $\widetilde{\textbf{x}}$ of $F$. Moreover, such a solution $\widetilde{\textbf{x}}$ must be nonsingular, i.e. $\textrm{rank}(\textrm{Jac}(F)|_{\widetilde{\textbf{x}}}) = n$. If one succeeds with this method using a conjugate symmetric box, i.e. $B=\overline{B}$, the unique root $\widetilde{\textbf{x}}$ must be real, since non-real solutions come in conjugate pairs.

The techniques of \textit{interval arithmetic} conservatively bound the image of $B$ under $N_F$ within another box $B'$ from which one proves that $B' \subset \textrm{interior}(B)$. This can be done using floating-point arithmetic in a certified manner.  See \cite{JuliaCertification} for details about the implementation of this, now standard, certification method. 

\subsection{Applying Krawczyk's method} We observe that the  polynomial system $F=(f_1,\ldots,f_{14}) = \textbf{0}$ is not square: it involves $14$ equations in $15$ variables. We write its solution set over $\mathbb{C}$ as $X_{\mathbb{C}}\subset \mathbb{C}^{15}$.
We introduce an additional linear equation $f_{15} = v_{1,1}-\frac{61}{17} = 0$
to our system. Note that the constant $\frac{61}{17}$ is a rational approximation to $(V^*)_{1,1} = 3.587996$.

Since $F=(f_1,\ldots,f_{15}) \in \mathbb{Q}[v_{1,1},\ldots,v_{3,5}]$ is a square system, one may apply Krawczyk's method by finding a box containing $V^*$ and using interval arithmetic to prove it contracts under $N_F$.  
We first certified $V^*$ via Krawczyk's method using the \texttt{Macaulay2} package \texttt{NumericalCertification} \cite{M2Certification} and subsequently in \texttt{julia}. 
After defining  $\texttt{F}$ and our candidate solution $\texttt{VR}$, the following portion of the \texttt{Macaulay2} code given in the Appendix executes in about five minutes.
\begin{lstlisting}[language=M2]
i31 : time cer = certifyRealSolution(F, VR)
 -- used 307.116s (cpu); 65.7524s (thread); 0s (gc)
o31 = true
\end{lstlisting}
In the \texttt{julia} implementation, the bounding box is explicitly returned:
\begin{lstlisting}[language=julia]
julia> @time C = certify(F,[VR])
 72.813907 seconds (7.26 M allocations: 1.559 GiB, 0.50% gc time)
CertificationResult
===================
- 1 solution candidates given
- 1 certified solution intervals (1 real, 0 complex)
- 1 distinct certified solution intervals (1 real, 0 complex)
julia> C.certificates[1].I
15x1 Arblib.AcbMatrix:
 [3.58823529411765 +/- 4.15e-15] + [+/- 9.22e-24]im
   [0.006182356890 +/- 3.66e-13] + [+/- 8.32e-15]im
     [4.5975381751 +/- 1.98e-11] + [+/- 3.11e-13]im
     [0.6880871518 +/- 4.12e-11] + [+/- 2.20e-13]im
   [-1.45618086541 +/- 5.87e-12] + [+/- 1.34e-13]im
    [-0.3624667542 +/- 5.06e-11] + [+/- 1.42e-13]im
    [1.87083079377 +/- 6.49e-12] + [+/- 1.42e-13]im
   [-0.37760827121 +/- 2.18e-12] + [+/- 4.88e-14]im
   [-0.22423435991 +/- 4.72e-12] + [+/- 6.62e-14]im
  [-0.130229923447 +/- 9.59e-13] + [+/- 1.93e-14]im
   [-4.26517137975 +/- 2.68e-12] + [+/- 6.36e-14]im
     [1.8104349519 +/- 5.58e-11] + [+/- 3.06e-13]im
    [-0.6782689763 +/- 1.49e-11] + [+/- 1.99e-13]im
    [-4.5378217507 +/- 3.48e-11] + [+/- 2.46e-13]im
   [-0.52932644272 +/- 5.79e-12] + [+/- 1.28e-13]im
\end{lstlisting}
This is a computational proof of the following theorem. 
\begin{theorem}
    \label{thm:compresult}
    There exists a real nonsingular point $\widetilde{V}$ on a one-dimensional real component of the solution set $X_{\mathbb{C}}$ of $F=\textbf{0}$. 
\end{theorem}

Theorem \ref{thm:compresult} implies that there exists a nonzero $3 \mathbb{Z} \times 5 \mathbb{Z}$-periodic potential which is Fermi isospectral to the zero potential at energy $\lambda = 0$ by Theorem~\ref{thm2}, proving Theorem \ref{thm1}. Since there exists a nontrivial  
$3\mathbb{Z}\times5\mathbb{Z}$-periodic real potential $\widetilde{V}$ such that
\begin{equation}\label{ga}
F_0(\widetilde{V})=F_0(\mathbf 0).
\end{equation}
Equation \eqref{gv} implies that
$[\widetilde{V}]=0$. Thus, $\widetilde{V}$ cannot be a constant function and  Theorem \ref{thm3} then follows from Remark  \ref{rem}  and equation
\eqref{ga}. Consequently, we answer Question~\ref{question1} in the negative and disprove Conjecture~\ref{conj1}.

\section*{Acknowledgments}
T. Brysiewicz  is supported by NSERC Discovery Grant RGPIN-2023-03551.  M. Faust was partially supported by NSF DMS-2052519. W. Liu was supported in part by NSF grants DMS-2246031 and DMS-2052572,
by a Simons Fellowship in Mathematics, and by a Visiting Miller Professorship
from the Miller Institute for Basic Research in Science, University of California, Berkeley.
Some computations used for this research were performed on the Texas A\&M University Department of Mathematics Whistler cluster, access to which was supported by NSF DMS-2201005.

 \bibliographystyle{abbrv}
\bibliography{absence}
\makeatletter
\@setaddresses
\global\let\@setaddresses\relax
\makeatother

\newpage
\section*{Appendix}
\noindent The following code justifies all computational claims made in this article.

\begin{center}
    \textbf{Macaulay2 - generate polynomials and certify}
\end{center}
\begin{lstlisting}[language=M2]
--------------------------------
------GENERATE EQUATIONS--------
--------------------------------
R = QQ[x,y,xi,yi,v_(1,1)..v_(3,5)] 
--(x,y,xi,yi) play the roles of (z_1,z_2,z_1^(-1),z_2^(-1)) resp.

DV = matrix {{v_(1,1), 1, xi, 1, 0, 0, 0, 0, 0, 0, 0, 0, y, 0, 0}, 
             {1, v_(2,1), 1, 0, 1, 0, 0, 0, 0, 0, 0, 0, 0, y, 0},
             {x, 1, v_(3,1), 0, 0, 1, 0, 0, 0, 0, 0, 0, 0, 0, y},
             {1, 0, 0, v_(1,2), 1, xi, 1, 0, 0, 0, 0, 0, 0, 0, 0},
             {0, 1, 0, 1, v_(2,2), 1, 0, 1, 0, 0, 0, 0, 0, 0, 0}, 
             {0, 0, 1, x, 1, v_(3,2), 0, 0, 1, 0, 0, 0, 0, 0, 0}, 
             {0, 0, 0, 1, 0, 0, v_(1,3), 1, xi, 1, 0, 0, 0, 0, 0}, 
             {0, 0, 0, 0, 1, 0, 1, v_(2,3), 1, 0, 1, 0, 0, 0, 0}, 
             {0, 0, 0, 0, 0, 1, x, 1, v_(3,3), 0, 0, 1, 0, 0, 0}, 
             {0, 0, 0, 0, 0, 0, 1, 0, 0, v_(1,4), 1, xi, 1, 0, 0}, 
             {0, 0, 0, 0, 0, 0, 0, 1, 0, 1, v_(2,4), 1, 0, 1, 0}, 
             {0, 0, 0, 0, 0, 0, 0, 0, 1, x, 1, v_(3,4), 0, 0, 1}, 
             {yi, 0, 0, 0, 0, 0, 0, 0, 0, 1, 0, 0, v_(1,5), 1, xi}, 
             {0, yi, 0, 0, 0, 0, 0, 0, 0, 0, 1, 0, 1, v_(2,5), 1}, 
             {0, 0, yi, 0, 0, 0, 0, 0, 0, 0, 0, 1, x, 1, v_(3,5)}};
                 
ZERO = join({x,y,xi,yi},new List from 15:0); 			
--substitutions for x,y,xi,yi,v_(1,1)..v_(3,5) to obtain D_0(x,y)
D0 = sub(DV,for i from 0 to numgens(R)-1 list R_i=>ZERO#i) 	
--apply substitutions
F0 = det(D0) 							
--take determinant
FV = det(DV, Strategy=>Cofactor); 				
--take determinant with no substitutions using cofactor expansion
FDifferenceModInverses = (F0-FV) % ideal(x*xi-1, y*yi-1); 	
--compute the difference modulo inverse relations
Eqs = unique(entries(
         (coefficients(FDifferenceModInverses,
                        Variables=>{x,xi,y,yi}))_1_0
      )); 							
--extract the coefficients of the difference 
--------------------------------------
----CERTIFY PROPOSED REAL SOLUTION----
--------------------------------------
needsPackage "NumericalCertification";

precBits = 200; RRring = RR_precBits;
T = RRring[for i from 4 to numgens(R)-1 list R_i]
F = for e in Eqs list sub(e,T) 
--substitute equations into the correct ring (no x,xi,y,yi)

VR = {3.58799648386156,  0.688112385771553,  1.87080104041807,
    -0.130206885382426, -0.678204748466961,  0.00618746135457902,
     -1.45598519094981, -0.377624670108577, -4.26499086966595,
     -4.53749591961956,   4.59717237464526, -0.36252210536849,
    -0.224238038666354,   1.81035923653508, -0.529360554357988
     };
F = polySystem (F | {v_(1,1)-61/17}); 

CCo = coefficientRing ring F;  VRCC = apply(VR, t -> t_CCo);
M = matrix{VRCC}; P0 = point M;
--apply some iterations of Newton
P1 = newton(F, P0); P2 = newton(F, P1); P3 = newton(F, P2);
VR = matrix P3;

time cer = certifyRealSolution(F, VR)

--------------------------------
----WRITE RESULTS TO A FILE-----
--------------------------------
--Macaulay2 Format
f = openOut "35M2Polys.txt";
for eq in Eqs do (f << toString eq << endl;);
close f;
--julia Format
f = openOut "35JuliaPolys.txt";
scan(Eqs, eq -> (
    s := toString eq;
    for i from 1 to 3 do for j from 1 to 5 do
        s = replace("v_\\(" | toString i | "," | toString j | "\\)",
                    "v[" | toString i | "," | toString j | "]", s);
    f << s << endl;));
close f;
\end{lstlisting}

\begin{center}
    \textbf{Julia - read polynomials, certify, and solve}
\end{center} 
\begin{lstlisting}[language=julia]
# Run in a julia session in same directory as 35JuliaPolys.txt
using HomotopyContinuation
@var v[1:3,1:5]; @var a[1:15]; @var c; V = vcat(v...);

############################################
########### Read and Certify ###############
############################################
S = [eval(Meta.parse(r)) for r in readlines("35JuliaPolys.txt")];
VR = [3.58799648386156,   0.00618746135457902, 4.59717237464526,
      0.688112385771553, -1.45598519094981,   -0.362522105368490,
      1.87080104041807,  -0.377624670108577,  -0.224238038666354,
     -0.130206885382426, -4.26499086966595,    1.81035923653508,
     -0.678204748466961, -4.53749591961956,   -0.529360554357988];

F = System(vcat(S,[v[1,1]-61//17]));
@time C = certify(F,[VR])

############################################
######## Solve for a real solution #########
############################################
I = solve(F; start_system = :total_degree, iterator_only=true)

global i=0; global total_successes = 0; 
function is_real_success(s)
    global i; i+=1; global total_successes;
    r=is_real(s); t=is_success(s); println(i,"  ",r,",",t)
    if t
        total_successes+=1; proportion = total_successes/i;
        println("Total success proportion:", proportion)
        println("Expected solution count:",4572288000*proportion)
    end
    return(r && t)
end
#The above function will print, for each path tracked, whether it
# was (a) successful and (b) real. It will estimate the expected
# proportion of start paths which are successful based on this data

J = Iterators.filter(s->is_real_success(s), I);
@time first(J)


\end{lstlisting}
\end{document}